\documentclass[12pt]{article}
\usepackage{amsfonts}
\usepackage{amssymb}
\usepackage{amsfonts,amssymb, amsmath, latexsym}
\usepackage{pst-all}
\usepackage{CJK}
\usepackage{mathrsfs}

\usepackage[left,pagewise,mathlines]{lineno}

\newtheorem{cor}{Corollary}[section]
\newtheorem{lem}{Lemma}[section]
\newtheorem{thm}{Theorem}[section]

\newtheorem{obs}[thm]{Observation}

\newtheorem{fact}{Fact}
\newtheorem{remark}{Remark}

\newenvironment{pf}[1][Proof]{\noindent\textbf{#1.} }{\hfill\rule{1mm}{2mm}}

\def\theequation{\thesection.\arabic{equation}}
\makeatletter \@addtoreset{equation}{section} \makeatother

\renewcommand\linenumberfont{\normalfont\bfseries\scriptsize}

\def\f{\phi}
\def\g{\gamma}
\def\sq{\mathbin{\square}}
\def\sse{\subseteq}
\def\spe{\supseteq}
\def\a{\alpha}
\def\b{\beta}
\def\br{\bar}
\def\th{\theta}
\def\tm{\times}
\def\vn{\varnothing}
\def\half{\frac 12}
\def\r{\rho}
\def\sm{\setminus}
\def\m{\mathcal}
\def\ma{\mathscr}

\begin{document}

\title{Trees with Maximum $p$-Reinforcement Number\thanks{The work was
supported by NNSF of China (No.10711233) and the Fundamental Research Fund of NPU (No. JC201150)}}
\author{
{ You Lu$^a$
\quad Jun-Ming Xu$^b$} \\
{\small $^a$Department of Applied Mathematics,}\\
          {\small     Northwestern Polytechnical University,}\\
 {\small             Xi'an Shaanxi 710072, P. R. China}\\
 {\small  Email: luyou@nwpu.edu.cn}\\ \\
{\small $^b$Department of Mathematics,}\\
{\small  University of Science and Technology of China,}   \\
{\small  Wentsun Wu Key Laboratory of CAS}\\
{\small Hefei, Anhui, 230026, P. R. China}\\
{\small Email: xujm@ustc.edu.cn}  \\
}

\date{}

\maketitle

\begin{abstract}
Let $G=(V,E)$ be a graph and $p$ a positive integer. The $p$-domination number
$\g_p(G)$ is the minimum cardinality of a set $D\subseteq V$ with $|N_G(x)\cap D|\geq p$ for all $x\in V\setminus D$. The
$p$-reinforcement number $r_p(G)$ is the smallest number of edges
whose addition to $G$ results in a graph $G'$ with $\g_p(G')<
\g_p(G)$. Recently, it was proved by Lu et al. that $r_p(T)\leq p+1$ for a tree $T$ and $p\geq 2$.
In this paper, we characterize all trees
attaining this upper bound for $p\geq 3$.

\end{abstract}




\section{Induction}

For notation and graph-theoretical terminology not defined here we
follow \cite{cl96,x03}. Let $G=(V,E)$ be a finite, undirected and simple
graph with vertex-set $V=V(G)$ and edge-set $E=E(G)$. For a vertex $x\in V$,
its {\it open neighborhood}, {\it closed
neighborhood} and {\it degree} are respectively
$N_{G}(x)=\{y\in V : xy\in E\}$, $N_{G}[x]=N_{G}(x)\cup \{x\}$ and
$deg_G(x)=|N_G(x)|$. A
vertex of degree one is called an {\it endvertex} or a {\it leaf} and its neighbor is called a {\it stem}. We denote the set of leaves of $G$ by $L(G)$.

For $S\subseteq V(G)$, the
subgraph induced by $S$ (resp. $V(G)\setminus S$) is denoted by $G[S]$
(resp. $G-S$). The complement $G^c$ of $G$ is the simple graph whose vertex-set is $V(G)$ and whose edges are the pairs of nonadjacent vertices of $G$. For $B\subseteq E(G^c)$, we use
$G+B$ to denote the subgraph with vertex-set $V(G)$ and edge-set $E(G)\cup
B$. For convenience, for vertex $v\in V(G)$, edge $xy\in E(G^c)$ and subgraph $H\subseteq G$, we write
$G-\{v\}$, $G+\{xy\}$ and $G-V(H)$ for $G-v$, $G+xy$ and $G-H$, respectively.

Let $T$ be a tree and $xy\in E(T)$.  Use the notation $T_y$ to denote
the component of $T-x$ containing $y$. To simplify notation, we still
use $T_y$ to denote $V(T_y)$. If $T$ is
a rooted tree, then, for every $x\in V(T)$, we let $C(x)$ and $D(x)$ denote the sets of
children and descendants, respectively, of $x$, and define $D[x] = D(x)\cup \{x\}$.

Let $p$ be a positive
integer. A subset $D\subseteq V(G)$ is a {\it $p$-dominating set} of
$G$ if every vertex not in $D$ has at least $p$ neighbors in $D$. The {\it
$p$-domination number} $\gamma_p(G)$ is the minimum cardinality of a
$p$-dominating set of $G$. A $p$-dominating set with
cardinality $\gamma_p(G)$ is called a $\gamma_p(G)$-set.
The {\it $p$-reinforcement number} $r_p(G)$ is
the smallest number of edges of $G^c$
that have to be added to $G$ in order to
reduce $\gamma_p(G)$, that is
 $$
 r_p(G)=\min\{|B|\ :\ B\subseteq E(G^c) \mbox{\ \ with\ \ } \g_p(G+B)< \g_p(G)\}.
 $$
By convention $r_p(G)=0$ if $\g_p(G)\leq p$. Clearly, the $1$-domination and $1$-reinforcement numbers are the well-known domination and reinforcement numbers, respectively.

The $p$-domination was introduced by Fink and Jacobson \cite{fj85} and has been well studied in graph theory (see, for example, \cite{bcf05,bcv06,cr90,dghp,f85,fhv08,hmv11}). Very recently, Chellali et al. \cite{cfhv11} gave an excellent survey on this topics. The $p$-reinforcement number introduced by Lu et al.~\cite{lhx12} is a parameter
for measuring the vulnerability of the $p$-domination, is also a natural extension of the classical reinforcement number which was introduced by Kok and Mynhardt \cite{km90} and has been studied by a number of authors
including \cite{cs03,hx12,hrr11,hwx09,zls03}. Motivated by the works of these authors, Lu et al.~\cite{lhx12} give an original study on the $p$-reinforcement for any $p\geq 1$. For a graph $G$ and $p\geq 1$, they found a method to determine $r_p(G)$ in terms of $\g_p(G)$, show that the
decision problem on $r_p(G)$ is NP-hard and established some upper bounds.

To be surprising, for a tree $T$, the upper bounds of $r_p(T)$ have distinct difference between $p=1$ and $p\geq 2$. For $p=1$, Blair et al. \cite{bgh08}
proved that $r_1(T)\leq \frac{1}{2}|V(T)|$ and this bound is sharp.
However, the following theorem implies that no result in terms of $|V(T)|$ exists for $p\geq 2$.

\begin{thm}\label{thm1.2}\textnormal{(\cite{lhx12})}\ \
$r_p(T)\leq p+1$ for a tree $T$ and $p\geq 2$.
\end{thm}

In this paper, we continue to consider the $p$-reinforcement number of trees. We will focus on the structural properties of the extremal trees in Theorem \ref{thm1.2}, and character all extremal trees for $p\geq 3$ by a recursive construction.

The rest of this paper is organized as follows. In Section 2 we present some notations and known results. We show the structural properties of a tree $T$ with $r_p(T)=p+1$ for $p\geq 3$ in Sections 3, and then
characterize such trees in Section 4.


\section{Known Results}
In this section, we will make the necessary preparations for proving the main results in Sections
3 and 4. Let $G=(V,E)$ be a graph and $p$ be a positive integer.

\begin{obs}\label{obs1.2}
Every $p$-dominating set contains all vertices of degree at most $p-1$.
\end{obs}

 For $X\subseteq
V$, let
 $
 X^*=\{x\in V\setminus X: |N_G(x)\cap X|< p\},
 $
and define
 \begin{eqnarray}
 \eta_p(x,X,G)&=&\left\{  \begin{array}{ll}
                        p-|N_G(x)\cap X| & \mbox{if }x\in X^*;\\
                        0 & \mbox{otherwise},
                      \end{array}
                  \right. \mbox{ for $x\in V$,}\label{e2.1}\\
 \eta_p(S, X, G)&=&\sum_{x\in S}\eta_p(x,X,G) \mbox{\ \ for $S\subseteq V$, }\label{e2.2} \\
 \eta_p(G)&=&\min\{\eta_p(V,X,G) : |X|<\g_p(G)\}\label{e2.3}.
 \end{eqnarray}

A subset $X\subseteq V$ is called an \emph{$\eta_p(G)$-set} if
$\eta_p(G)=\eta_p(V,X,G)$.

\begin{obs}\label{obs2.1}
If $X$ is an $\eta_p(G)$-set, then $|X|=\g_p(G)-1$.
\end{obs}


\begin{thm}\label{thm2.2}\textnormal{(\cite{lhx12})}\ \
For a graph $G$ and $p\geq 1$, $r_p(G)=\eta_p(G)$ if
$\g_p(G)\geq p+1$.
\end{thm}

\begin{cor}\label{cor2.3}\textnormal{(\cite{lhx12})}\ \
Let $G$ be a graph and $p\geq 1$. For any $H\subseteq G$ with $\g_p(H)\geq p+1$ and $\g_p(G)\geq \g_p(H)+\g_p(G-H)$,
$$r_p(G)\leq r_p(H).$$
\end{cor}

 Let $X\subseteq V$ and $x\in
X$. A vertex $y\in V\setminus X$ is called a {\it $p$-private
neighbor} of $x$ with respect to $X$ if $xy\in E$ and
$|N_G(y)\cap X|=p$.  The {\it $p$-private
neighborhood} of $x$ with respect to $X$, denoted by $PN_p(x,X,G)$, is defined as
the set of $p$-private neighbors of $x$ with respect to $X$.
Let
 \begin{eqnarray}
 \mu_p(x,X,G)&=&|PN_p(x,X,G)|+\max\{0,p-|N_G(x)\cap X|\},\label{e2.4}\\
 \mu_p(X,G)&=&\min\{\mu_p(x,X,G) : x\in X\},              \label{e2.5}\\
 \mu_p(G)&=&\min\{\mu_p(X,G):\ X \mbox{ is a $\g_p(G)$-set}\}\label{e2.6}.
\end{eqnarray}

\begin{thm}\label{thm2.4}\textnormal{(\cite{lhx12})}\ \
For a graph $G$ and $p\geq 1$,
 $
 r_p(G)\leq \mu_p(G)
 $
with equality if $r_p(G)=1$.
\end{thm}

\section{Structural Theorems}

In this section, we investigate some structural properties of
trees with $r_p(T)=p+1$. We need the following lemma which is
given in \cite{lh11}.

\begin{lem}\label{lem3.1}\textnormal{(\cite{lh11})}\
Let $p\geq 2$ and $T$ be a tree with a $p$-dominating set $D$. Then
$D$ is a unique $\g_p(T)$-set if and only if $D$ satisfies either
$|N_T(x)\cap D|\leq p-2$ or $|PN_p(x,D,T)|\geq 2$ for each $x\in D$
with degree at least $p$.
\end{lem}

\begin{thm}\label{thm3.2}\
Let $p\geq 2$ and $T$ be a tree with $r_p(T)=p+1$.  For a $\g_p(T)$-set $D$,
\begin{description}
  \item[\ \  \ (i)] $PN_p(x,D,T)\ne\emptyset$ for each $x\in D$;
  \item[\ \  (ii)] $D$ is a unique $\g_p(T)$-set.
\end{description}
\end{thm}

\begin{pf}\ Let $x$ be a vertex in $D$.
By the definitions of $\mu_p$ in (\ref{e2.4})$\sim$(\ref{e2.6})
and Theorem \ref{thm2.4}, we have
\begin{eqnarray*}
|PN_p(x,D,T)|&=&\mu_p(x,D,T)-\max\{0,p-|N_T(x)\cap D|\}\\
&\geq &\mu_p(D,T)-p\\
&\geq& \mu_p(T)-p\\
&\geq& r_p(T)-p\\
&=&1.
\end{eqnarray*}
The first conclusion holds.

We now prove the second conclusion. It is obvious that $|D|=\g_p(T)\geq p+1$ since $r_p(T)=p+1>0$. If every vertex in $D$ has degree less than $p$,
then $D$ is a unique $\g_p(T)$-set by Observation~\ref{obs1.2}.
Assume that $D$ has a vertex, say $x$, with degree at least $p$. By the definitions of
$\mu_p$ in (\ref{e2.4})$\sim$(\ref{e2.6}) and Theorem \ref{thm2.4}, we have
 \begin{eqnarray}\label{e3.0}
 r_p(T)\leq\mu_p(T)\leq\mu_p(x,D,T)= |PN_p(x,D,T)|+\max\{0,p-|N_T(x)\cap D|\}.
 \end{eqnarray}
If $|N_T(x)\cap D|\geq p-1$, then $\max\{0,p-|N_T(x)\cap D|\}\leq 1$, and so we obtain from (\ref{e3.0}) that
 $
 |PN_p(x,D,T)|\geq r_p(T)-1=p\geq 2.
 $
This fact implies that $D$ satisfies the conditions of Lemma \ref{lem3.1}. Thus $D$ is a unique $\g_p(T)$-set. The theorem follows.
\end{pf}

Through this paper, for any tree $T$ with a unique $\g_p(T)$-set, we denote the unique $\g_p(T)$-set by $\ma{D}_T$ for short.

\begin{thm}\label{thm3.3}
Let $p\geq 2$ and $T$ be a tree with $r_p(T)=p+1$. For edge $xy\in E(T)$ with $x\in
\ma{D}_T$, let $T_y$ denote the component of $T-x$ containing $y$. Then
\begin{description}
  \item[\ \ \ (i)] If $y\in PN_p(x,\ma{D}_T,T)$, then either $T_y$ is a star $K_{1,p-1}$
  with center $y$, or $r_p(T_y)=1$ and $\ma{D}_T\cap T_y$ is an $\eta_p(T_y)$-set.
  Moreover, $\eta_p(T_y,S,T_y)\geq p-1$ for $S\subseteq T_y$ with $|S|=|\ma{D}_T\cap T_y|$ and $y\in S$.
\item[\ \ (ii)] If $y\notin PN_p(x,\ma{D}_T,T)$, then $r_p(T_y)=p+1$
  and $\ma{D}_{T_y}=\ma{D}_T\cap T_y$.
\end{description}
\end{thm}

\begin{pf}
\textbf{(i)} Since $x\in \ma{D}_T$ and $y\in
PN_p(x,\ma{D}_T,T)$, we have $y\notin \ma{D}_T$
and $|N_{T}(y)\cap \ma{D}_T|=p$.
It follows that
$$|T_y|\geq |((N_{T}(y)\cap \ma{D}_T)\setminus \{x\})\cup \{y\}|=p.$$
If $|T_y|=p$, then it is easy to see that $T_y=K_{1,p-1}$ with center $y$,
and $\eta_p(T_y,S,T_y)=p-1$ for any $S\subseteq T_y$ with
$|S|=|\ma{D}_T\cap T_y|$ and $y\in S$.
Assume that $|T_y|\geq p+1$ in the following.

\begin{claim}\label{claim1}
 $\g_p(T_y)\geq p+1$.
\end{claim}

\noindent\textbf{Proof of Claim \ref{claim1}.} Suppose, to be contrary, that $\g_p(T_y)=p$, and let $D_y$ be a $\g_p(T_y)$-set. Then there is a vertex, say $z$, of $T_y$ not in $D_y$ since $|T_y|\geq p+1$. To $p$-dominate $z$, $|N_{T_y}(z)\cap D_y|\geq p=|D_y|\geq 2$. From the fact that two vertices in a tree have at most one common neighbor,  we know that $z$ is a unique vertex of $T_y$ not in $D_y$. Thus $T_y=K_{1,p}$ and $z=y$. Hence, by Observation \ref{obs1.2}, we can obtain a contradiction as follows:
$$p=|N_T(y)\cap \ma{D}_T|=|N_T(y)|=|(T_y\setminus \{y\})\cup \{x\}|=1+p.$$
The claim holds. $\Box$

\begin{claim} \label{claim2}
$
\eta_p(T_y,S,T_y)\geq \left\{\begin{array}{ll}
                                  1 & \mbox{if\ \ $y\notin S$};\\
                                  p-1 & \mbox{if\ \ $y\in S$},
                              \end{array}
                      \right.
$
for any $S\subseteq T_y$ with $|S|=|\ma{D}_T\cap T_y|$.
\end{claim}

\noindent \textbf{Proof of Claim \ref{claim2}.} Let $S$ be a subset of $T_y$ such that $|S|=|\ma{D}_T\cap T_y|$. We complete the proof of Lemma \ref{claim2} by distinguishing the following three cases.

If $S=\ma{D}_T\cap T_y$, then it is obvious that $\eta_p(T_y,S,T_y)=1$ by $y\in PN_p(x,\ma{D}_T,T)$.

If $S\neq \ma{D}_T\cap T_y$ and $y\notin S$, then let $D'=(\ma{D}_T\setminus T_y)\cup S$. Clearly, $|D'|=|\ma{D}_T|=\g_p(T)$ and $D'\neq \ma{D}_T$. Since $\ma{D}_T$ is a unique $\g_p(T)$-set, $D'$ is not a $p$-dominating set of $T$. Thus, by $x\in \ma{D}_T$  and $x\in D'$, we have
$$\eta_p(T_y,S,T_y)\geq \eta_p(T_y,D',T)=\eta_p(V(T),D',T)\geq 1.$$

If $y\in S$, then, for any $u\in N_T(y)$, we use $T_u$ to denote the component of $T-y$ containing $u$.  Since $y\in
S\setminus \ma{D}_T$ and
$|S|=|\ma{D}_T\cap T_y|$, we have
  \begin{eqnarray*}
 \sum_{u\in N_T(y)\setminus \{x\}}|S\cap T_u|= |S\setminus \{y\}|&=&|S|-1\\
  &=&|\ma{D}_T\cap T_y|-1
  =\sum_{u\in N_T(y)\setminus \{x\}}|\ma{D}_T\cap T_u|-1.
  \end{eqnarray*}
Thus, there is some $v\in N_T(y)\setminus \{x\}$ such that $|S\cap T_v|\leq |\ma{D}_T\cap T_v|-1$. Let
 $$
 D''=(\ma{D}_T\setminus T_v)\cup (S\cap T_v).
 $$
Then
 \begin{eqnarray*}
 |D''|&\leq&|\ma{D}_T\setminus T_v|+|T_v\cap S|\\
 &\leq& |\ma{D}_T\setminus T_v|+(|\ma{D}_T\cap T_v|-1)\\
 &=&|\ma{D}_T|-1\\
 &=&\g_p(T)-1,
  \end{eqnarray*}
  and $\eta_p(y,D'',T)\leq 1$ since $|N_T(y)\cap D''|\geq|N_T(y)\cap \ma{D}_T-\{v\}|\geq |N_T(y)\cap \ma{D}_T|-1=p-1$.
By Theorem \ref{thm2.2} and the definitions of $\eta_p$ in (\ref{e2.1})$\sim$(\ref{e2.3}), we have
\begin{eqnarray*}
    p+1=r_p(T)=\eta_p(T)&\leq& \eta_p(V(T),D'',T)\\
                      &=&\eta_p(T_v,D'',T)+\eta_p(y,D'',T)\\
                      &\leq &\eta_p(T_v,D'',T)+1\\
                      &=&\eta_p(T_v,T_v\cap S,T_y)+1.
 \end{eqnarray*}
Thus, $\eta_p(T_v, S\cap T_v, T_y)\geq p$, and so
  \begin{eqnarray*}
 \eta_p(T_y,S,T_y)&\geq& \eta_p(T_v,S,T_y)\\
  &=& \eta_p(T_v, S\cap T_v,T_y)-1\\
  &\geq& p-1.
  \end{eqnarray*}
 The proof of Claim \ref{claim2} is complete. $\Box$

Claim \ref{claim2} implies that the second conclusion in \textbf{(i)} holds. We now show that $r_p(T_y)=1$ and $\ma{D}_T\cap T_y$ is an $\eta_p(T_y)$-set.
By $p\geq 2$, it is easily seen from  Claim \ref{claim2} that $\eta_p(T_y,S,T_y)\geq 1$ for $S\subseteq T_y$ with $|S|=|\ma{D}_T\cap T_y|$. So $\g_p(T_y)\geq |\ma{D}_T\cap T_y|+1$, further, we have $\g_p(T_y)=|\ma{D}_T\cap T_y|+1$ since $(\ma{D}_T\cap T_y)\cup \{y\}$ is a $p$-dominating set of $T_y$. By Claim \ref{claim2} and its proof, $\eta_p(T_y)=\eta_p(T_y,\ma{D}_p\cap T_y,T_y)=1$ and $\ma{D}_T\cap T_y$ is an $\eta_p(T_y)$-set.  Hence $r_p(T_y)=\eta_p(T_y)=1$ by Theorem \ref{thm2.2} and Claim \ref{claim1}.

The proof of \textbf{(i)} is complete.

\textbf{(ii)}  By the hypothesis of $y\notin
PN_p(x,\ma{D}_T,T)$, $\ma{D}_T\setminus T_y$ and
$\ma{D}_T\cap T_y$ are $p$-dominating sets of $T-T_y$ and
$T_y$, respectively. Thus, we have
 \begin{equation}\label{e3.1}
 |\ma{D}_T\setminus T_y|\geq \g_p(T-T_y)\mbox{\ \ and\ \ }|\ma{D}_T\cap T_y|\geq \g_p(T_y).
 \end{equation}
Note that the union of a $\g_p(T-T_y)$-set and a $\g_p(T_y)$-set is a $p$-dominating set of $T$. Thus, we have
 \begin{equation}\label{e3.2}
 \g_p(T-T_y)+\g_p(T_y)= \g_p(T).
 \end{equation}
It follows from (\ref{e3.1}) and (\ref{e3.2}) that
 \begin{eqnarray*}
  \g_p(T-T_y)+\g_p(T_y)=\g_p(T)=|\ma{D}_T|&=&|\ma{D}_T\setminus T_y|+|\ma{D}_T\cap T_y|\\
  &\geq& \g_p(T-T_y)+\g_p(T_y),
 \end{eqnarray*}
which yields that $|\ma{D}_T\cap T_y|=\g_p(T_y)$, that is,
$\ma{D}_T\cap T_y$ is a $\g_p(T_y)$-set.

To the end, by \textbf{(ii)} of
Theorem~\ref{thm3.2}, it is sufficient to prove $r_p(T_y)=p+1$.

Arbitrary take $u\in \ma{D}_T\cap T_y$. Since $r_p(T)=p+1$, we have
$PN_p(u,\ma{D}_T,T)\ne \emptyset$ by \textbf{(i)} of Theorem \ref{thm3.2}. Let $z\in PN_p(u,\ma{D}_T,T)$.
Clearly, $z\in T_y$ and $z\ne y$ since $u\in T_y$ and $y\notin PN_p(x,\ma{D}_T,T)$. So
 $$
 |N_{T_y}(z)\cap (\ma{D}_T\cap T_y)|=|N_T(z)\cap \ma{D}_T|=p.
 $$
It follows that
 \begin{equation*}
 \g_p(T_y)=|\ma{D}_T\cap T_y|\geq |N_{T_y}(z)\cap (\ma{D}_T\cap T_y)|=p.
 \end{equation*}

Furthermore, we can show that $\g_p(T_y)\geq p+1$.
To be contrary, assume that $\g_p(T_y)=p$. Then $|\ma{D}_T\cap T_y|=p$ and  $z$ is a unique vertex of $T_y$ not in $\ma{D}_T\cap T_y$ since two vertices in a tree have at most one common neighbor. It follows that $T_y=K_{1,p}$ with center $z$ and $\ma{D}_T\cap T_y=T_y\setminus \{z\}$, and so $y\in \ma{D}_T$ and $y$ is a leaf of $T_y$. This implies that $N_T(y)=\{x,z\}$ and by (\ref{e2.4}),
 \begin{eqnarray*}
 \mu_p(y,\ma{D}_T,T)&=&|PN_p(y,\ma{D}_T,T)|+\max\{0,p-|N_T(y)\cap
 \ma{D}_T|\}\\
 &=& 1+(p-1)=p.
 \end{eqnarray*}
By Theorem \ref{thm2.4} and the definitions of $\mu_p$ in
(\ref{e2.4})$\sim$(\ref{e2.6}), we have
 $$
 r_p(T)\leq \mu_p(T)\leq \mu_p(y,\ma{D}_T,T)=p,
 $$
which contradicts to the hypothesis of $r_p(T)=p+1$.
Thus $\g_p(T_y)\geq p+1$.

Let $X$ be an $\eta_p(T_y)$-set and $Y=X\cup (\ma{D}_T\setminus T_y)$. Then $|X|=\g_p(T_y)-1$ and
 $$
 |Y|=(\g_p(T_y)-1)+(\g_p(T)-\g_p(T_y))= \g_p(T)-1.
 $$
By Theorem \ref{thm2.2}, we have
 \begin{eqnarray*} r_p(T_y)&=&\eta_p(T_y,X,T_y)\\
                         &\geq&\eta_p(T_y,Y,T)= \eta_p(V(T),Y,T)\\
                         &\geq& \eta_p(T)= r_p(T)=p+1.
 \end{eqnarray*}
Combining this with Theorem~\ref{thm1.2}, we have $r_p(T_y)=p+1$,
and so \textbf{(ii)} is true.

The theorem follows.\end{pf}

\begin{lem}\label{lem3.4}
Assume that $p\geq 3$ and $T$ be a tree with $r_p(T)=p+1$. Let $x\in \ma{D}_T$  and $X\subseteq V(T-x)$ with $|X|<\g_p(T)$. If $\mu_p(x,\ma{D}_T,T)\geq p+2$ and $\eta_p(V(T),X,T)=p+1$, then
$|X\cap T_y|= |\ma{D}_T\cap T_y|$ for all $y\in N_T(x)$, where $T_y$ is the component of $T-x$ containing $y$.
\end{lem}

\begin{pf}
\ Let $N_T(x)=\{x_1,\ldots,x_d\}$, where $d=deg_T(x)$, and let $T_i$ be the component of
$T-x$ containing $x_i$. Combining $x\in \ma{D}_T\setminus X$ with $|X|<\g_p(T)$, we have
$$\sum_{j=1}^d|X\cap T_j|=|X|\leq \g_p(G)-1=|\ma{D}_T\setminus \{x\}|=\sum_{j=1}^d|\ma{D}_T\cap T_j|.$$
We only need to show that $|X\cap T_j|=|\ma{D}_T\cap T_j|$ for any $j\in\{1,\ldots,d\}$.
Suppose, to the contrary, that there exists some
$i\in \{1,\ldots,d\}$ such that
 \begin{equation*}
 |X\cap T_i|< |\ma{D}_T\cap T_i|.
 \end{equation*}
Our aim is to deduce a contradiction.

We first give the following two claims. Claim \ref{claim3} can be easily obtained from Theorem \ref{thm3.3},.

\begin{claim}\label{claim3}
For any $j\in \{1,\ldots,d\}$,
 $$
 |\ma{D}_T\cap T_j|= \left\{\begin{array}{ll}
                           \g_p(T_j)-1 & \mbox{if\ \ $x_j\in PN_p(x,\ma{D}_T,T)$};\\
                           \g_p(T_j)   &  \mbox{if\ \ $x_j\notin PN_p(x,\ma{D}_T,T)$}.
                           \end{array}
                  \right.
 $$
\end{claim}

\begin{claim}\label{claim4}
For any $j\in \{1,\ldots,d\}$,
  $$
 |X\cap T_j|\geq \left\{\begin{array}{ll}
                           \g_p(T_j)-2 & \mbox{if\ \ $j=i$};\\
                           \g_p(T_j)   &  \mbox{if\ \ $j\ne i$}.
                           \end{array}
                  \right.
 $$

\end{claim}

\noindent \textbf{Proof of Claim \ref{claim4}.} Let
 $
 D=(X\cap T_i)\cup (\ma{D}_T\setminus T_i).
 $
Clearly,
 \begin{eqnarray*}
 |D|=|X\cap T_i|+|\ma{D}_T\setminus T_i|
 <|\ma{D}_T\cap T_i|+|\ma{D}_T\setminus T_i|
 =|\ma{D}_T|=\g_p(T).
 \end{eqnarray*}
Thus $\eta_p(V(T),D,T)\geq \eta_p(T)$ by (\ref{e2.3}). From $x\in \ma{D}_T\setminus T_i\subseteq D$, we know that $D$ is a $p$-dominating set of $T-T_i$, and so $\eta_p(T_i,D,T)=\eta_p(V(T),D,T)$. Note that $\g_p(T)\geq p+1$ by $r_p(T)= p+1$. Together with the hypothesis of Lemma \ref{lem3.4}, (\ref{e2.1})$\sim$(\ref{e2.3}) and Theorem \ref{thm2.2}, we have
 \begin{eqnarray}
 p+1&=&\eta_p(V(T),X,T)\nonumber \\
 &\geq& \eta_p(V(T),X,T)-\sum\limits_{j\neq i}\eta_p(T_j,X,T)-\eta_p(x,X,T)\nonumber\\
 &=&\eta_p(T_i,X,T)\label{e3.3}\\
 &\geq& \eta_p(T_i,D,T)=\eta_p(V(T),D,T)\nonumber\\
 &\geq& \eta_p(T)=r_p(T)\nonumber\\
 &=&p+1,\nonumber
 \end{eqnarray}
which implies that all the equalities in (\ref{e3.3}) hold. In
particular,
 $$
 \mbox{$\eta_p(x,X,T)=\eta_p(T_j,X,T)=0$,\ \ for $j\neq i$,}
 $$
which means that
 \begin{eqnarray}
 & & |N_T(x)\cap X|\geq p\ \ {\rm and}\label{e3.4}\\
 & & |N_{T_j}(u)\cap X|\geq p,\ \ {\rm for}\  j\neq i\ {\rm and}\ u\in T_j\setminus X.
 \label{e3.5}
 \end{eqnarray}
It follows from (\ref{e3.5}) and $x\notin X$ that $X\cap T_j$
(for $j\neq i$) is a $p$-dominating set of $T_j$, and so
 \begin{equation*}
 |X\cap T_j|\geq \g_p(T_j)\ \ {\rm for}\ j\neq i.
 \end{equation*}

Also since all the equalities in (\ref{e3.3}) hold, we have
$\eta_p(V(T),D,T)=r_p(T)$, which means that $D$ is an $\eta_p(T)$-set. By Observation \ref{obs2.1},
 \begin{equation}\label{e3.6}
 |D|=\g_p(T)-1=|\ma{D}_T|-1.
 \end{equation}
Since $x\in \ma{D}_T$ and $x_i\in N_T(x)$, $(\ma{D}_T\cap
T_i)\cup\{x_i\}$ is a $p$-dominating set of $T_i$, and so
 \begin{equation}\label{e3.7}
 |\ma{D}_T\cap T_i|\geq \g_p(T_i)-1.
 \end{equation}
It follows from the definition of $D$, (\ref{e3.6}) and
(\ref{e3.7}) that
 \begin{eqnarray*}
 |X\cap T_i|&=&|D|-|\ma{D}_T\setminus T_i|\\
                   &=&(|\ma{D}_T|-1)-|\ma{D}_T\setminus T_i|\\
                   &=&|\ma{D}_T\setminus (\ma{D}_T\setminus T_i)|-1\\
                   &=&|\ma{D}_T\cap T_i|-1\\
                   &\geq& \g_p(T_i)-2.
                   \end{eqnarray*}
Claim \ref{claim4} holds. $\Box$

Together with Claims \ref{claim4} and \ref{claim3}, we have that
 \begin{eqnarray}
\g_p(T)&\geq&|X|+1=\sum\limits_{j\neq i}|X\cap T_j|+|X\cap T_i|+1\nonumber\\
           &\geq& \sum\limits_{j\neq i}\g_p(T_j)+(\g_p(T_i)-2)+1\nonumber\\
           &=&\sum\limits_{i=1}^d\g_p(T_j)-1\nonumber\\
           &=&\sum\limits_{x_j\in PN_p(x,\ma{D}_T,T)}\g_p(T_j)+\sum\limits_{x_j\notin PN_p(x,\ma{D}_T,T)}\g_p(T_j)-1\label{e3.8}\\
           &=&\sum\limits_{x_j\in PN_p(x,\ma{D}_T,T)}(|\ma{D}_T\cap T_j|+1)+\sum\limits_{x_j\notin PN_p(x,\ma{D}_T,T)}(|\ma{D}_T\cap T_j|)-1 \nonumber\\
           &=&|PN_p(x,\ma{D}_T,T)|+\sum\limits_{i=1}^d|\ma{D}_T\cap T_j|-1\nonumber\\
           &=& |PN_p(x,\ma{D}_T,T)|+(|\ma{D}_T|-1)-1\ \ \ \ ({\rm since}\ x\in \ma{D}_T)\nonumber\\
           &=& |PN_p(x,\ma{D}_T,T)|+\g_p(T)-2,\nonumber
\end{eqnarray}
which yields
 \begin{equation}\label{e3.9}
 |PN_p(x,\ma{D}_T,T)|\leq 2.
 \end{equation}
It follows from our hypothesis of $\mu_p(x,\ma{D}_T,T)\geq p+2$
and (\ref{e3.9}) that
 \begin{eqnarray*}
 p+2&\leq&\mu_p(x,\ma{D}_T,T)\\
 &=&|PN_p(x,\ma{D}_T,T)|+\max\{0,p-|N_T(x)\cap \ma{D}_T|\}\\
 &\leq& 2+\max\{0,p-|N_T(x)\cap \ma{D}_T|\}\\
 &\leq& p+2,
 \end{eqnarray*}
which implies that
\begin{eqnarray}
&&|N_T(x)\cap \ma{D}_T|=0 \mbox{\ \ and}\label{e3.10}\\
  &&|PN_p(x,\ma{D}_T,T)|=2.\label{e3.11}
 \end{eqnarray}
It follows from (\ref{e3.11}) that all equalities in (\ref{e3.8}) hold, and so, by Claim \ref{claim4},
  \begin{equation}\label{e3.12}
  |X\cap T_j|=\left\{
 \begin{array}{ll}
 \g_p(T_i)-2 &\ \mbox{ if $j=i$};\\
 \g_p(T_j) &\ \mbox{ if $j\neq i$}.
 \end{array}
 \right.
 \end{equation}

We claim that $x_i\in PN_p(x,\ma{D}_T,T)$. Assume, to be contrary, that $x\notin PN_p(x,\ma{D}_T,T)$, then $r_p(T_i)=p+1$ by \textbf{(ii)} of Theorem \ref{thm3.3}. Since $\eta_p(T_i,X\cap T_i,T_i)=\eta_p(T_i,X,T)=p+1$ by $x\notin X$ and (\ref{e3.3}),
$X\cap T_i$ is an $\eta_p(T_i)$-set, and so $|X\cap T_i|=\g_p(T_i)-1$ by Observation \ref{obs2.1}, which contradicts with (\ref{e3.12}). The claim holds, and so if $x_j\notin PN_p(x,\ma{D}_T,T)$ then $j\ne i$.

For $j\ne i$, it can be easily seen from (\ref{e3.12}) that $X\cap T_j$ is a $\g_p(T_j)$-set  since $X\cap T_j$ is a $p$-dominating set of $T_j$ by (\ref{e3.5}) and $x\notin X$. On the other hand, for any $x_j\notin PN_p(x,\ma{D}_T,T)$, $\ma{D}_T\cap T_j$ is a unique $\g_p(T_j)$-set by \textbf{(ii)} of Theorem~\ref{thm3.3}. Hence
 \begin{equation}\label{e3.13}
 X\cap T_j=\ma{D}_T\cap T_j,\ \ \mbox{for $x_j\notin PN_p(x,\ma{D}_T,T)$}.
 \end{equation}

Together with (\ref{e3.10}), (\ref{e3.11}) and (\ref{e3.13}), we have
$|N_T(x)\cap X|\leq |PN_p(x,\ma{D}_T,T)|=2,$
which contradicts to (\ref{e3.4}) since $p\geq 3$. The lemma follows.
\end{pf}

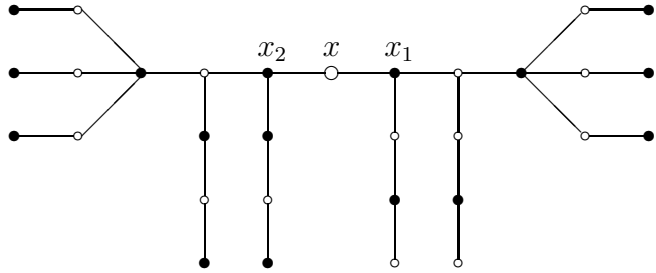
\begin{figure}[h]
\begin{center}
\psset{unit=0.8pt}

\begin{pspicture}(400,200)(0,40)

\put(50,150){\circle*{5}}\put(80,150){\circle{4}}\put(110,150){\circle*{5}}
\put(140,150){\circle{4}}\put(170,150){\circle*{5}}\put(200,150){\circle{6}}
\put(230,150){\circle*{5}}\put(260,150){\circle{4}}\put(290,150){\circle*{5}}
\put(320,150){\circle{4}}\put(350,150){\circle*{5}}

\put(50,180){\circle*{5}}\put(50,120){\circle*{5}}\put(80,180){\circle{4}}\put(80,120){\circle{4}}
\put(320,180){\circle{4}}\put(350,180){\circle*{5}}\put(320,120){\circle{4}}\put(350,120){\circle*{5}}

\put(140,120){\circle*{5}}\put(140,90){\circle{4}}\put(140,60){\circle*{5}}
\put(170,120){\circle*{5}}\put(170,90){\circle{4}}\put(170,60){\circle*{5}}
\put(230,120){\circle{4}}\put(230,90){\circle*{5}}\put(230,60){\circle{4}}
\put(260,120){\circle{4}}\put(260,90){\circle*{5}}\put(260,60){\circle{4}}

\put(196,158){$x$}
\put(165,158){$x_2$}
\put(225,158){$x_1$}

\put(52,150){\line(1,0){26}}\put(82,150){\line(1,0){26}}\put(112,150){\line(1,0){26}}
\put(142,150){\line(1,0){26}}\put(172,150){\line(1,0){25}}\put(203,150){\line(1,0){25}}
\put(232,150){\line(1,0){26}}\put(262,150){\line(1,0){26}}\put(292,150){\line(1,0){26}}
\put(322,150){\line(1,0){26}}

\put(52,180){\line(1,0){26}}\put(52,120){\line(1,0){26}}
\put(82,180){\line(1,-1){28}}\put(82,120){\line(1,1){28}}
\put(322,180){\line(1,0){26}}\put(322,120){\line(1,0){26}}
\put(290,150){\line(1,1){28}}\put(290,150){\line(1,-1){28}}

\put(140,148){\line(0,-1){26}}\put(140,120){\line(0,-1){28}}
\put(140,88){\line(0,-1){26}}

\put(170,148){\line(0,-1){26}}\put(170,120){\line(0,-1){28}}
\put(170,88){\line(0,-1){26}}

\put(230,148){\line(0,-1){26}}\put(230,118){\line(0,-1){26}}
\put(230,88){\line(0,-1){26}}

\put(260,148){\line(0,-1){26}}\put(260,118){\line(0,-1){28}}
\put(260,88){\line(0,-1){26}}

 \end{pspicture}
\caption{\label{f1} \footnotesize A trees $T$ with  $\g_2(T)=17$ and $r_2(T)=3$. The vertex $x$ satisfies $\mu_2(x,\ma{D}_T,T)=4$ and the set $X$ of black vertices in $T$ has $|X|=16< \g_2(T)$ and $\eta_2(V(T),X,T)=3$. However, $|X\cap T_1|=7\ne 8=|\ma{D}_T\cap T_1|$ and $|X\cap T_2|=9\ne 8=|\ma{D}_T\cap T_2|$, where $T_1$ and $T_2$ are two components of $T-x$ containing $x_1$ and $x_2$, respectively.}
\end{center}
\end{figure}

\begin{remark}
The conclusion of Lemma \ref{lem3.4} may not be valid for $p=2$. An counterexample is showed in Figure \ref{f1}.
\end{remark}

\begin{thm}\label{thm3.5}
Let $p\geq 3$, $T$ be a tree with $r_p(T)=p+1$ and $x\in
\ma{D}_T$. If $\mu_p(x,\ma{D}_T,T)\geq p+2$, then
$\eta_p(V(T),X,T)\geq p+2$ for any $X\subseteq V(T-x)$ with
$|X|<\g_p(T)$.
\end{thm}

\begin{pf}\ Suppose, to the contrary, that $\eta_p(V(T),X,T)\leq p+1$ for some subset $X\subseteq V(T-x)$ with
$|X|<\g_p(T)$.
Since $\eta_p(V(T),X,T)\geq r_p(T)=p+1$ by Theorem~\ref{thm2.2}, we have
\begin{equation}\label{e3.15'}
\eta_p(V(T),X,T)= p+1.
\end{equation}

 Let
$T_y$ denote the component of $T-x$ containing $y$ for any $y\in N_T(x)$.
We partition $N_T(x)$ into six subsets $N_i$
($1\leq i\leq 6$) such that
 $N_T(x)=\bigcup_{i=1}^6N_i$,
 where
 $$
 \left\{\begin{array}{ll}
  N_1=PN_p(x,\ma{D}_T,T)\cap X,&
  N_2=PN_p(x,\ma{D}_T,T)- X, \\
  N_3=(N_T(x)\cap X)\cap \ma{D}_T,&
  N_4=(N_T(x)- X)\cap \ma{D}_T,\\
  N_5=(N_T(x)\cap X)-N_1-N_3, &
  N_6=(N_T(x)- X)-N_2-N_4.
  \end{array}\right\}
  $$
Clearly,
 \begin{equation}\label{e3.14}
 |N_1|+|N_2|=|PN_p(x,\ma{D}_T,T)|
 \end{equation}
 and, by $x\notin X$ and the definition of $\eta_p$ in (\ref{e2.1}),
 \begin{equation}\label{e3.15}
 \eta_p(x,X,T)=
 \max\{0,p-|N_1|-|N_3|-|N_5|\}.
 \end{equation}

For any $y\in N_1\cup N_2\cup N_5$, we can obtain $|X\cap T_y|=|\ma{D}_T\cap T_y|$ from (\ref{e3.15'}) and Lemma \ref{lem3.4}. If $y\in N_1\cup N_2$, then $y\in PN_p(x,\ma{D}_T,T)$ and, by \textbf{(i)} of Theorem \ref{thm3.3}, we have
 \begin{equation}\label{e3.16}
 \eta_p(T_y,X,T)=\eta_p(T_y,X\cap T_y,T_y)\geq \left\{\begin{array}{ll}
                                                                               p-1 & \mbox{if $y\in N_1$};\\
                                                                               1   &  \mbox{if $y\in N_2$}.
                                                                               \end{array}
                                                                               \right.
 \end{equation}
If $y\in N_5$, then $y\notin PN_p(x,\ma{D}_T,T)$ and $y\in X\setminus \ma{D}_T$, which implies that $X\cap T_y\ne \ma{D}_T\cap T_y$. By $|X\cap T_y|=|\ma{D}_T\cap T_y|$ and Theorem \ref{thm3.3} \textbf{(ii)}, $X\cap T_y$ is not a $p$-dominating set of $T_y$, and so
 \begin{equation}\label{e3.17}
 \eta_p(T_y,X,T)=\eta_p(T_y,X\cap T_y,T_y)\geq 1.
 \end{equation}

Therefore, we can derive from (\ref{e3.14})$\sim$(\ref{e3.17}) that
\begin{eqnarray*}
 \eta_p(V(T),X,T) &=& \eta_p(x,X,T)+\sum_{y\in N_T(x)}\eta_p(T_y,X,T) \\
   &\geq&  \max\{0,p-|N_1|-|N_3|-|N_5|\}+\sum_{y\in N_1\cup N_2\cup N_5}\eta_p(T_y,X,T)\\
   &\geq& \max\{0,p-|N_3|\}-|N_1|-|N_5|+(p-1)|N_1|+|N_2|+|N_5|\\
   &=& |PN_p(x,\ma{D}_T,T)|+\max\{0,p-|N_3|\}+(p-3)|N_1|\\
   &\geq& |PN_p(x,\ma{D}_T,T)|+\max\{0,p-|N_T(x)\cap\ma{D}_T|\}+0\\
   &=&\mu_p(x,\ma{D}_T,T)\\
   &\geq& p+2,
\end{eqnarray*}
which contradicts with (\ref{e3.15'}).
\end{pf}


\section{A constructive Characterization}

In this section, we will give a constructive characterization of
trees with $r_p(T)=p+1$ for $p\geq 3$.  To the end, assume that
$p\geq 3$.

\begin{figure}[h]
\begin{center}
\psset{unit=0.8pt}

\begin{pspicture}(400,180)(0,10)

\put(30,100){\circle*{6}}\put(25,85){$x$}
\put(80,125){\circle{5}}\put(70,133){$x_1$}
\put(80,75){\circle{5}}\put(70,63){$x_2$}

\put(30,100){\line(2,-1){48}}
\put(82,127){\line(2,1){28}}\put(82,123){\line(2,-1){28}}
\put(110,110){\circle*{5}}\put(110,140){\circle*{5}}\put(108,120){$\vdots$}

\put(30,100){\line(2,1){48}}
\put(82,77){\line(2,1){28}}\put(82,73){\line(2,-1){28}}
\put(110,60){\circle*{5}}\put(110,90){\circle*{5}}\put(108,70){$\vdots$}

\put(60,25){$F_{p-1}$}

\put(240,100){\circle{6}}\put(235,85){$x$}
\put(290,125){\circle*{5}}
\put(330,133){$x_1$}\put(340,125){\circle{5}}
\put(290,75){\circle*{5}}
\put(330,63){$x_t$}\put(340,75){\circle{5}}
\put(338,95){$\vdots$}

\put(242,98){\line(2,-1){47}}\put(290,125){\line(1,0){47}}
\put(342,127){\line(2,1){28}}\put(342,123){\line(2,-1){28}}
\put(370,110){\circle*{5}}\put(370,140){\circle*{5}}\put(368,120){$\vdots$}

\put(242,102){\line(2,1){47}}\put(290,75){\line(1,0){48}}
\put(342,77){\line(2,1){28}}\put(342,73){\line(2,-1){28}}
\put(370,60){\circle*{5}}\put(370,90){\circle*{5}}\put(368,70){$\vdots$}

\put(300,25){$F_{t,p-1}$}

 \end{pspicture}
\caption{\label{f2} \footnotesize Trees $F_{p-1}$ and $F_{t,p-1}$, where $t\geq p$ and each $x_i$ has $p-1$ leaves.}
\end{center}
\end{figure}
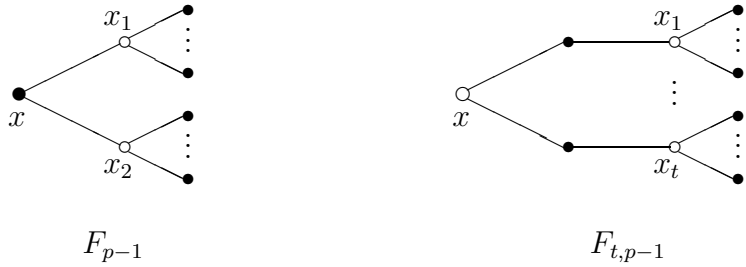

Let $t\geq 2$ be an integer. The unique stem of a star $K_{1,t}$ is called the {\it center} of $K_{1,t}$. A {\it spider} $S_t$ is a tree obtained by attaching one leaf at each endvertex of $K_{1,t}$.
Two important trees $F_{p-1}$ and $F_{t,p-1}$ in our construction are shown in Figure \ref{f2}, where $F_{p-1}$ (resp. $F_{t,p-1}$) is obtained by attaching $p-1$ leaves at  every endvertex of a path $P_3$ (resp. a spider $S_t$).

In Figure \ref{f2}, we call $x$ the \emph{center} of
$F_{p-1}$ and $F_{t,p-1}$.
Clearly, the set of black vertices in $F_{p-1}$ (resp. $F_{t,p-1}$) is a unique
$\g_p(F_{p-1})$-set (resp. $\g_p(F_{t,p-1})$-set). By Theorem
\ref{thm2.2}, the following lemma is  straightforward by computing
$\eta_p$.

\begin{lem}\label{lem4.1}
If $p\geq 3$, then
$r_p(F_{p-1})=p+1$ and $r_p(F_{t,p-1})=p+1$ for any $t\geq p$.
\end{lem}

Let $G$ and $H$ be two graphs. For $x\in V(G)$ and $y\in V(H)$, the notation $G\uplus_{xy} H$ is a graph consisting of $G$ and $H$ with an extra edge $xy$.
Let $T_0=K_{1,p}$ be a star and $A(T_0)=L(K_{1,p})$. We define a family $\mathscr{T}_p$ of trees as follows:
\begin{eqnarray*}
\mathscr{T}_p=\{T : T \mbox{ is obtained from a sequence $T_0, T_1, \ldots, T_k\ (k \geq 1)$ of trees, where\ \ }\\
\mbox{$T=T_k$ and $T_{i+1}\ (0\leq i\leq k-1)$ is obtained from $T_i$\ \ }\\
\mbox{by one of the operations listed below}\}.
\end{eqnarray*}
Four operations on a tree $T_i$:
\begin{flushleft}
\begin{itemize}
 \item[$\mathscr{O}_1$:]\mbox{$T_{i+1}=K_{1,p-1}\uplus_{xy}T_i$, where $x$ is the center of $K_{1,p-1}$ and $y\in A(T_i)$.}\\
                        \mbox{Let $A(T_{i+1})=L(K_{1,p-1})\cup A(T_i)$.}\\
 \item[$\mathscr{O}_2$:]\mbox{$T_{i+1}=K_{1,p}\uplus_{xy}T_i$, where $x$ is the center of $K_{1,p}$ and $y\notin A(T_i)$.}\\
                        \mbox{Let $A(T_{i+1})=L(K_{1,p})\cup A(T_i)$.}\\
 \item[$\mathscr{O}_3$:]\mbox{$T_{i+1}=F_{p-1}\uplus_{xy} T_i$, where $x$ is the center of $F_{p-1}$ and $y\in A(T_i)$ satisfying}\\
                        \mbox{$|PN_p(y,A(T_i),T_i)|\geq \min\{p+1,|N_{T_i}(y)\cap A(T_i)|+2$\}.}\\
                        \mbox{Let $A(T_{i+1})=\ma{D}_{F_{p-1}}\cup A(T_i)$.}\\
 \item[$\mathscr{O}_4$:]\mbox{$T_{i+1}=F_{t, p-1}\uplus_{xy} T_i$, where $t\geq p$, $x$ is the center of $F_{t, p-1}$ and $y\in T_i$.}\\
                        \mbox{Let $A(T_{i+1})=\ma{D}_{F_{t,p-1}}\cup A(T_i)$.}
\end{itemize}
\end{flushleft}

\begin{lem}\label{lem4.2}
Let $p\geq 3$ be an integer. If $T\in \mathscr{T}_p$, then $A(T)$ is a unique
$\g_p(T)$-set and $r_p(T)=p+1$.
\end{lem}

\begin{pf}\
Since $T\in \mathscr{T}_p$, there is an integer $k\geq 1$ such that
$T=T_k$ obtained from
a sequence $T_0,T_1, \ldots, T_{k-1},T_k$ of trees, where $T_0 = K_{1,p}$ and $T_{i+1}$ ($0\leq i\leq k-1$) can be obtained from $T_i$ by some $\mathscr{O}_j$ ($j\in \{1,2,3,4\}$). We complete the proof of Lemma \ref{lem4.2} by
induction on $k$.

For $k=1$, it is trivial that $A(T_1)$ is a
unique $\g_p(T_1)$-set and $r_p(T_1)=\eta_p(T_1)=p+1$ by straightly computing
$\eta_p(T_1)$. This establishes the induction base.

Let $k\geq 2$. Assume that
 $A(T_i)$ is a unique $\g_p(T_i)$-set and $r_p(T_i)=p+1$ for
$1\leq i\leq k-1$.
 Basing on this assumption, we will prove that $A(T_k)$ is a unique $\g_p(T_k)$-set and $r_p(T_k)=p+1$.

 Since $T_k$ is obtained from $T_{k-1}$ by $\mathscr{O}_j$ for some $j$, we denote $G_j=T_k-V(T_{k-1})$ and so $T_k=G_j\uplus_{xy} T_{k-1}$, where $x$ is the center of $G_j$ and $y$
satisfies the condition of $\mathscr{O}_j$.

We first show that $A(T_k)$ is a unique $\g_p(T_k)$-set. Applying the
induction on $T_{k-1}$, $A(T_{k-1})$ is a unique $\g_p(T_{k-1})$-set. For $v\in A(T_{k-1})$ with degree at least $p$, it follows from Lemma \ref{lem3.1} that either $|N_{T_{k-1}}(x)\cap A(T_{k-1})|\leq p-2$ or $|PN_p(v,A(T_{k-1}),T_{k-1})|\geq 2$. Hence, we can easily derive from the definition of $A(T_k)$ that $A(T_k)$ is a $p$-dominating set of $T_k$, moreover, satisfies the conditions of Lemma \ref{lem3.1}. By Lemma \ref{lem3.1}, $A(T_k)$ is a unique $\g_p(T_k)$-set

 Next we only need to prove
that $r_p(T_k)\geq p+1$ by Theorem~\ref{thm1.2}. Assume, to be
contrary, that $r_p(T_k)\leq p$. We will deduce a contradiction. Let $X$ be an $\eta_p(T_k)$-set such that $|X\cap V(T_{k-1})|$ is as large as possible. Then $|X\cap V(T_{k-1})|\leq \g_p(T_{k-1})$ by the choice of $X$ and, by Theorem \ref{thm2.2},
\begin{equation}\label{e4.1}
\eta_p(V(T_k),X,T_k)=r_p(T_k)\leq p.
\end{equation}
In addition, we also have the following three facts.

\begin{fact}\label{fact1}
At most one leaf of $T_k$ is not in $X$.
\end{fact}

Otherwise,
$\eta_p(V(T_k),X,T_k)\geq 2(p-1)> p$ by (\ref{e2.1}) and $p\geq 3$, which contradicts with (\ref{e4.1}).

\begin{fact}\label{fact2}
$|X\cap V(T_{k-1})|=\g_p(T_{k-1})-1$ or $\g_p(T_{k-1})$.
\end{fact}

Suppose, to be contrary, that $|X\cap V(T_{k-1})|\leq \g_p(T_{k-1})-2$. Since each $\eta_p(T_{k-1})$-set has cardinality $\g_p(T_{k-1})-1$ by Observation \ref{obs2.1}, $\eta_p(V(T_{k-1}), X\cap V(T_{k-1}),T_{k-1})\geq \eta_p(T_{k-1})+1$. Note the assumption that $r_p(T_{k-1})=p+1$. By (\ref{e2.1})$\sim$(\ref{e2.3}) and Theorem \ref{thm2.2}, we have
\begin{eqnarray*}
\eta_p(V(T_k),X,T_k)&\geq&\eta_p(V(T_{k-1}),X,T_k)\\
                    &\geq& \eta_p(V(T_{k-1}), X\cap V(T_{k-1}),T_{k-1})-1\\
                    &\geq& (\eta_p(T_{k-1})+1)-1\\
                    &=& r_p(T_{k-1})\\
                    &=&p+1.
\end{eqnarray*}
This contradicts with (\ref{e4.1}). Hence $|X\cap V(T_{k-1})|=\g_p(T_{k-1})-1$ or $\g_p(T_{k-1})$.

\begin{fact}\label{fact3}
If $|X\cap V(T_{k-1})|=\g_p(T_{k-1})-1$, then
$\{x,y\}\cap X=\{x\}$ and $\eta_p(V(G_j),X,T_k)=0$.
\end{fact}

Suppose, to be contrary, that $\{x,y\}\cap X\neq \{x\}$ or $\eta_p(V(G_j),X,T_k)> 0$.
From (\ref{e2.1})$\sim$(\ref{e2.3}) and Theorem \ref{thm2.2}, we will deduce a contradiction with (\ref{e4.1}).
If $\{x,y\}\cap X\neq \{x\}$, then
 \begin{eqnarray*}
 \eta_p(V(T_k),X,T_k)&\geq& \eta_p(V(T_{k-1}),X,T_k)\\
                     &=&\eta_p(V(T_{k-1}),V(T_{k-1})\cap X,T_{k-1})\\
                     &\geq& \eta_p(T_{k-1})=r_p(T_{k-1})=p+1,
 \end{eqnarray*}
 a contradiction with (\ref{e4.1}). If $\eta_p(V(G_j),X,T_k)> 0$, then
\begin{eqnarray*}
 \eta_p(V(T_k),X,T_k)&=&\eta_p(V(G_j),X,T_k)+\eta_p(V(T_{k-1}),X,T_k)\\
                     &>&\eta_p(V(T_{k-1}),X,T_k)\\
                     &\geq&\eta_p(V(T_{k-1}),X\cap V(T_{k-1}),T_{k-1})-1\\
                     &\geq& r_p(T_{k-1})-1\\
                     &=&p,
 \end{eqnarray*}
 which is also a contradiction with (\ref{e4.1}). Hence Fact \ref{fact3} holds.

Based on the above facts, we complete the proof of Lemma \ref{lem4.2} by distinguishing the following two cases.

\textbf{Case 1.} $j=1$ or $2$.

Then $G_j=K_{1,\ell}$ with center $x$, where $\ell=p-1$ if $j=1$, otherwise $\ell=p$. Note that $X$ is an $\eta_p(T_k)$-set. By Observation \ref{obs2.1}, we have
\begin{equation}\label{e4.2}
|X|=\g_p(T_k)-1=|A(T_k)|-1=|A(T_{k-1})|+\ell-1=\g_p(T_{k-1})+\ell-1.
\end{equation}
By Fact \ref{fact2}, we only need to consider the following two subcases.

If $|X\cap V(T_{k-1})|=\g_p(T_{k-1})$, then $|X\cap V(G_j)|=\ell-1$ by (\ref{e4.2}), that is, exactly two vertices are not in $X$. Since $G_j=K_{1,\ell}$ is a star with center $x$, it follows from Fact \ref{fact1}
that $x\notin X$ and there is a unique leaf $x'$ of $G_j$
not in $X$. Thus
\begin{eqnarray*}
p  &=& \eta_p(x',X,T_k)\\
  &=& \eta_p(V(T_k),X,T_k)-\eta_p(V(T_{k-1}),X,T_k)-\eta_p(x,X,T_k)\\
  &\leq& p-\eta_p(V(T_{k-1}),X,T_k)-\eta_p(x,X,T_k),
\end{eqnarray*}
which implies that $\eta_p(V(T_{k-1}),X,T_k)=\eta_p(x,X,T_k)=0$.
Together with $x'\notin X$, $x\notin X$ and $\eta_p(x,X,T_k)=0$, we have
$$
|N_{T_k}(x)\setminus \{x'\}|\geq |N_{T_k}(x)\cap X|\geq p,
$$
and so $\ell=p$ (this means that $j=2$ and $T_k$ is obtained from $T_{k-1}$ by $\mathscr{O}_2$)
and $y\in X$.
By $y\in X$,
 $$
 \eta_p(V(T_{k-1}),X\cap V(T_{k-1}),T_{k-1})=\eta_p(V(T_{k-1}),X,T_k)=0,
 $$
which implies that $X\cap V(T_{k-1})$ is a $p$-dominating set of $T_{k-1}$. Further, $X\cap V(T_{k-1})=A(T_{k-1})$ since $|X\cap V(T_{k-1})|=\g_p(T_{k-1})$ and $A(T_{k-1})$ is a unique $\g_p(T_{k-1})$-set. By the condition of $\ma{O}_2$, we have $y\notin A(T_{k-1})=X\cap V(T_{k-1})$, a contradiction with $y\in X$.

If $|X\cap V(T_{k-1})|=\g_p(T_{k-1})-1$, then $|X\cap V(G_j)|=\ell$ by (\ref{e4.2})
and $x\in X$ by Fact \ref{fact3}. Since $G_j=K_{1,\ell}$ is a star with center $x$, there is a leaf of $G_j$
not in $X$. Hence $\eta_p(V(G_j),X,T_k)\geq p-1>0$, which contradicts
with the result $\eta_p(V(G_j),X,T_k)=0$ in Fact \ref{fact3}.

\textbf{Case 2.} $j=3$ or $4$.

Then $G=F_{p-1}$ or $G_{t,p-1}$ ($t\geq p$) with center $x$. Since $X$ is an $\eta_p(T_k)$-set, it is easily seen from Observation \ref{obs2.1} that \begin{equation}\label{e4.3}
|X|=\g_p(T_k)-1=|A(T_k)|-1=\g_p(G_j)+|A(T_{k-1})|-1=\g_p(G_j)+\g_p(T_{k-1})-1.
\end{equation}

\textbf{Subcase 2.1.} $|X\cap V(T_{k-1})|=\g_p(T_{k-1})$.

It is clear that $|X\cap V(G_j)|=\g_p(G_j)-1$ by (\ref{e4.3}). Since $r_p(G_j)=p+1$ by Lemma \ref{lem4.1},
we can deduce $\{x,y\}\cap X=\{y\}$ and $\eta_p(V(T_{k-1}),X,T_k)=0$
by the similar proof of Fact \ref{fact3}.

 If $j=3$, then $G_j=F_{p-1}$ with center $x$ and $\mu_p(x,\ma{D}_{G_j},G_j)=p+2$ by (\ref{e2.4}). Applying Theorem \ref{thm3.5}, we obtain that $\eta_p(V(G_j),X\cap V(G_j),G_j)\geq p+2$ from $|X\cap V(G_j)|<\g_p(G_j)$ and $x\notin X$. By (\ref{e4.1}), we get a contradiction that
\begin{eqnarray*}
 p&\geq& \eta_p(V(T_k),X,T_k)=\eta_p(V(G_j),X,T_k)\\
 &\geq& \eta_p(V(G_j),X\cap V(G_j),G_j)-1\\
 &\geq& p+1.
 \end{eqnarray*}

If $j=4$, then $G_j=F_{t,p-1}$ ($t\geq p$) with center $x$ and, by Figure \ref{f2}, every vertex of $\ma{D}_{G_j}-\{x\}$ has degree $1$ in $G_j-x$. Thus there
is at most one vertex in $\ma{D}_{G_j}$ but not in $X$
(otherwise, we have $\eta_p(V(T_k),X,T_k)=\eta_p(V(G_j),X,T_k)\geq
2(p-1)>p$ by $x\notin X$ and $p\geq 3$, which contradicts with (\ref{e4.1})). Combining this with $|X\cap V(G_j)|=\g_p(G_j)-1=|\ma{D}_{G_j}|-1$, we have
 $X\cap V(G_j)\subset \ma{D}_{G_j}$. Let $\ma{D}_{G_j}-X=\{u\}$ and $v$ be the unique neighbor of $u$ in $G_j-x$. Then
 \begin{eqnarray*}
 \eta_p(V(T_k),X,T_k)&\geq &\eta_p(u,X,T_k)+\eta_p(v,X,T_k)\\
                     &=&\eta_p(u,X\cap V(G_j),G_j)+\eta_p(v,X\cap V(G_j),G_j)\\
                     &\geq& p+1.
 \end{eqnarray*}
This contradicts with (\ref{e4.1}).

\textbf{Subcase 2.2.} $|X\cap V(T_{k-1})|=\g_p(T_{k-1})-1$.

 By Fact \ref{fact3}, $\{x,y\}\cap X=\{x\}$ and $\eta_p(V(G_j),X,T_k)=0$.
So $X\cap V(G_j)$ is a
$p$-dominating set of $G_j$, moreover, has cardinality $\g_p(G_j)$ by (\ref{e4.3}). This means that
$X\cap V(G_j)$ is a $\g_p(G_j)$-set. So $X\cap V(G_j)=\ma{D}_{G_j}$ since $\ma{D}_{G_j}$ is a unique
$\g_p(G_j)$-set. By $x\in X$ and Figure \ref{f2}, we have
$G_j=F_{p-1}$, that is, $T_k$ is obtained from $T_{k-1}$ by
$\mathscr{O}_3$. Thus $y\in A(T_{k-1})$
by the condition of $\ma{O}_3$.

We first show that there is a neighbor $z$ of $y$ in $A(T_{k-1})$ but not in $X$. For this aim, we consider that $|N_{T_{k-1}}(y)\cap A(T_{k-1})|$ and $|N_{T_{k-1}}(y)\cap X|$.

 Applying the inductive assumption on $T_{k-1}$, $A(T_{k-1})$ is a unique $\g_p(T_{k-1})$-set and $r_p(T_{k-1})=p+1$. It follows from (\ref{e2.4})$\sim$(\ref{e2.6}) and Theorem \ref{thm2.4} that $\mu_p(y,A(T_{k-1}),T_{k-1})\geq \mu_p(T_{k-1})\geq r_p(T_{k-1})= p+1$.
 Further, we can claim that
 \begin{equation}\label{e4.4}
 \mu_p(y,A(T_{k-1}),T_{k-1})= p+1.
 \end{equation}
 Assume, to the contrary, that
$\mu_p(y,A(T_{k-1}),T_{k-1})\geq p+2$. On the one hand, by $r_p(T_{k-1})=p+1$ and Theorem~\ref{thm3.5}, we have
$\eta_p(V(T_{k-1}),X\cap V(T_{k-1}),T_{k-1})\geq p+2$, and so
 $$
 \eta_p(V(T_{k-1}),X,T_k)\geq\eta_p(V(T_{k-1}),X\cap V(T_{k-1}),T_{k-1})-1\geq p+1.
 $$
 On the other hand, by (\ref{e4.1}) and $\eta_p(V(G_j),X,T_k)=0$, we have
 $$\eta_p(V(T_{k-1}),X,T_k)=\eta_p(V(T_k),X,T_k)-\eta_p(V(G_j),X,T_k)\leq p,$$
 a contradiction. Thus (\ref{e4.4}) holds.

 For convenience, let $A_y=N_{T_{k-1}}(y)\cap A(T_{k-1}).$  Together with (\ref{e4.4}), (\ref{e2.4}) and the condition of $\ma{O}_3$, we have
 \begin{eqnarray*}
 p+1&=&\mu_p(y,A(T_{k-1}),T_{k-1})\\
     &=&|PN_p(y,A(T_{k-1}),T_{k-1})|+\max\{0,p-|A_y|\}\\
     &\geq& \left\{\begin{array}{ll}
                          (p+1)+0             & \mbox{if\ \ $|A_y|\geq p$};\\
                          (|A_y|+2)+(p-|A_y|) & \mbox{if\ \ $|A_y|< p$}
                    \end{array}
             \right.\\
     &=& \left\{\begin{array}{ll}
                          p+1 & \mbox{if\ \ $|A_y|\geq p$};\\
                          p+2 & \mbox{if\ \ $|A_y|< p$},
                    \end{array}
             \right.
 \end{eqnarray*}
which implies that
\begin{equation}\label{e4.5}
|N_{T_{k-1}}(y)\cap A(T_{k-1})|=|A_y|\geq p \mbox{\ \ and\ \ } |PN_p(y,A(T_{k-1}),T_{k-1})|=p+1.
\end{equation}

Since $|X\cap V(T_{k-1})|=\g_p(T_{k-1})-1$ and $r_p(T_{k-1})=p+1$, it follows from (\ref{e2.1})$\sim$(\ref{e2.3}) and Theorem \ref{thm2.2} that
 $
 \eta_p(V(T_{k-1}),X\cap V(T_{k-1}),T_{k-1})\geq \eta_p(T_{k-1})=r_p(T_{k-1})= p+1.
 $
By (\ref{e4.1}), we have
 \begin{equation}\label{e4.6}
 |N_{T_{k-1}}(y)\cap X|< p.
 \end{equation}

 (\ref{e4.5}) and (\ref{e4.6}) implies that there is a neighbor $z$ of $y$ in
$A(T_{k-1})$ but not in $X$. Let $T_z$ be the component
of $T_{k-1}-y$ containing $z$. By $z\in A(T_{k-1})$, we have $z\notin
PN_p(y,A(T_{k-1}),T_{k-1})$. By
\textbf{(ii)} of Theorem~\ref{thm3.3}, $r_p(T_z)=p+1$ and $A(T_{k-1})\cap V(T_z)$ is
a unique $\g_p(T_z)$-set.

If $|X\cap V(T_z)|\geq \g_p(T_z)$, then let $X'=(X- V(T_z))\cup
(A(T_{k-1})\cap V(T_z))$. Clearly, $|X'|\leq|X|$. Together with
$y\notin X$, $z\in A(T_{k-1})-X$ and
(\ref{e4.6}), we have
 \begin{eqnarray*}
 \eta_p(V(T_{k}),X',T_{k})&=&\eta_p(V(T_k)- \{y\}- V(T_z),X',T_k)+\eta_p(y,X',T_k)+\eta_p(V(T_z),X',T_k)\\
  &=&\eta_p(V(T_k)- \{y\}- V(T_z),X,T_k)+(\eta_p(y,X,T_k)-1)+0\\
  &=&\eta_p(V(T_k)-V(T_z),X,T_k)-1\\
  &\leq&\eta_p(V(T_k),X,T_k)-1
 \end{eqnarray*}
  which contradicts that $X$ is an $\eta_p(T_{k-1})$-set.

If $|X\cap V(T_z)|<\g_p(T_z)$, then, by $y\notin X$ and $r_p(T_z)=p+1$, we have
\begin{eqnarray*}
\eta_p(V(T_k),X,T_{k})&\geq& \eta_p(V(T_z),X,T_k)\\
&=&\eta_p(V(T_z),X\cap V(T_z),T_z)\\
&\geq& \eta_p(T_z)= r_p(T_z)=p+1,
\end{eqnarray*}
which contradicts with (\ref{e4.1}).

This completes the proof of Lemma \ref{lem4.2}.
\end{pf}

Let $p\geq 1$ be an integer. For a tree $T$ with a unique $\g_p(T)$-set $\ma{D}_T$, we use $\ell_p(T)$ to
denote the number of all $p$-private neighbors with respect to
$\ma{D}_T$ in $T$. Since every $p$-private neighbor with respect to $\ma{D}_T$ has exactly $p$ neighbors in $\ma{D}_T$, we have
 \begin{equation}\label{e4.7}
 \ell_p(T)=\frac{1}{p}\sum_{x\in \ma{D}_T}|PN_p(x,\ma{D}_T,T)|.
 \end{equation}

\begin{lem}\label{lem4.3}
Let $p\geq 3$ be an integer and $T$ be a tree. If $r_p(T)=p+1$, then $T\in \mathscr{T}_p$.
\end{lem}

\begin{pf}
By Theorem~\ref{thm3.2}, $T$ has a unique $\g_p(T)$-set $\ma{D}_T$
and each vertex in $\ma{D}_T$ has at least one $p$-private
neighbor with respect to $\ma{D}_T$.
By $|\ma{D}_T|=\g_p(T)> p$ and (\ref{e4.7}), we have
 $$
 \ell_p(T)=\frac{1}{p}\sum_{x\in \ma{D}_T}|PN_p(x,\ma{D}_T,T)|\geq\frac{1}{p}\sum_{x\in \ma{D}_T}1= \frac{1}{p}|\ma{D}_T|> 1,
 $$
that is, $\ell_p(T)\geq 2$ since $\ell_p(T)$ is an integer. We will show $T\in \mathscr{T}_p$ by induction on $\ell_p(T)$.

If $\ell_p(T)=2$, then let $x$ and $y$ be two $p$-private neighbors with
respect to $\ma{D}_T$. Since every vertex in $\ma{D}_T$ has at least one $p$-privated neighbor with respect to $\ma{D}_T$,
 we have
$$\ma{D}_T=(N_T(x)\cap \ma{D}_T)\cup (N_T(y)\cap \ma{D}_T).$$
From the fact that $x$ and $y$ have at
most one common neighbor in $T$, we know that
 $$
 \g_p(T)=|\ma{D}_T|=\left\{\begin{array}{ll}
                              2p-1  &  \mbox{if\ \ $|N_T(x)\cap N_T(y)|=1$};\\
                              2p    &   \mbox{if\ \ $|N_T(x)\cap N_T(y)|=0$}.
                              \end{array}
                              \right.
 $$
Combining with $r_p(T)=p+1$, we can check easily that
$$
T=\left\{  \begin{array}{ll}
                        F_{p-1} & \mbox{if\ \ } \g_p(T)=2p-1; \\
                        S_{p,p} & \mbox{if\ \ } \g_p(T)=2p,
                      \end{array}
                  \right.
$$
where $S_{p,p}$ is obtained by attaching $p-1$ leaves
at every endvertex of a path $P_2$. Clearly, $T$ can be obtained from star $T_0=K_{1,p}$ by
$\mathscr{O}_1$ if $T=F_{p-1}$, otherwise by
$\mathscr{O}_2$. So $T\in \mathscr{T}_p$. This establishes the base
case.

Let $\ell_p(T)\geq 3$. Assume that, for any tree $T'$ with $r_p(T')=p+1$, if $\ell_p(T')< \ell_p(T)$ then $T'\in
\mathscr{T}_p$ .

For a leaf $r$ of $T$, we root $T$ at $r$.  Let $d=\max\{d_T(r,x) : x\in V(T)\}$. Thus we can partite $V(T)$ into $\{V_0,V_1,\ldots,V_d\}$, where
 $$
 V_i=\{x\in V(T)\ \ :\ \ \max\{d_T(x,v) : v\in D[x]\}=i\}, \mbox{\ \ \ for $i=0,1,\ldots, d$.}
 $$
By the definitions of $V_i$, $V_0=L(T)-\{r\}$ and $V_d=\{r\}$. Clearly, $V_0\subseteq \ma{D}_T$ by Observation \ref{obs1.2} and $d\geq 3$ (otherwise,
$T$ is a star, a contradiction with $\ell_p(T)\geq 3$).

We now turn to consider $V_1$. For any $x\in V_1$, we know from the definition of $V_1$ that $x$ is a stem of $T$ and $C(x)\subseteq V_0$. By $r_p(T)=p+1$, $x$ is a
$p$-private neighbor with
respect to $\ma{D}_T$, and so $deg_T(x)=p $ or $p+1$.

\begin{claim}\label{claim5}
If there exists a vertex $x\in V_1$ with $deg_T(x)=p+1$, then $T\in \mathscr{T}_p$.
\end{claim}

\noindent \textbf{Proof of Claim \ref{claim5}.} Let $T_x=T[D[x]]$ and
$T'=T-D[x]$. By $deg_T(x)=p+1$, $T_x=K_{1,p}$ with center $x$, and so $T=K_{1,p}\uplus_{xy}T'$, where $y$ is the father
of $x$ in $T$.

Since $|C(x)|=deg_T(x)-1=p$ and $x$ is a $p$-private neighbor with respect to $\ma{D}_T$, we
have $y\notin \ma{D}_T$ and $C(x)$ (resp. $\ma{D}_T\cap V(T')$) is a $p$-dominating set of $T_x$ (resp. $T'$). So
 \begin{equation}\label{e4.8}
 \g_p(T)=|\ma{D}_T|=|C(x)|+|\ma{D}_T\cap V(T')|\geq\g_p(T_x)+\g_p(T').
 \end{equation}
Further, we have $\g_p(T)= \g_p(T_x)+\g_p(T')$ since the union
between a $\g_p(T_x)$-set and a $\g_p(T')$-set is a $p$-dominating set of $T$. Hence, $\ma{D}_T\cap V(T')$ is a $\g_p(T')$-set.

Since $x\notin \ma{D}_T$ and $y\notin \ma{D}_T$, we have $\g_p(T')\geq p$, further, $\g_p(T')\geq p+1$ by $\ell_p(T)\geq 3$.
Together with Theorem \ref{thm1.2}, (\ref{e4.8}) and Corollary \ref{cor2.3}, we have
 $$p+1\geq r_p(T')\geq r_p(T)=p+1,$$
which implies that $r_p(T')=p+1$. By Theorem \ref{thm3.2},
$V(T')\cap \ma{D}_T$ is a unique $\g_p(T')$-set. Noting that $x$ is a $p$-private neighbor with respect to $\ma{D}_T$, we have $\ell_p(T')=\ell_p(T)-1$.
 Applying the
induction on $T'$, $T'\in \mathscr{T}_p$ and $\ma{D}_{T'}=\ma{D}_T\cap V(T')=A(T')$ by Lemma \ref{lem4.2}. Since $T=K_{1,p}\uplus_{xy}T'$ and $y\notin \ma{D}_T\cap V(T')=A(T')$, $T$ is obtained from $T'$ by
$\mathscr{O}_2$, and so $T\in \mathscr{T}_p$.
$\Box$

In the following, by Claim \ref{claim5}, we only need to consider the case that
\begin{equation}\label{e4.9}
deg_T(v)=p, \mbox{\ \ for each $v\in V_1$.}
\end{equation}
 Then the father of each vertex in $V_1$ belongs to
$\ma{D}_T$, and so $V_2\subseteq \ma{D}_T$.

Let $x\in V_3$
and $P=xwvu$ be a path in $T[D[x]]$ such that $deg_T(w)$ is as large
as possible. Obviously, $u\in V_0$, $v\in V_1$ and $w\in V_2$. By (\ref{e2.4})$\sim$(\ref{e2.6}) and
Theorem \ref{thm2.4},
 $$
 \mu_p(w,\ma{D}_T,T)\geq \mu_p(T)\geq r_p(T)=p+1.
$$

\textbf{Case 1.} $\mu_p(w,\ma{D}_T,T)\geq p+2$.

Let $T'=T-D[v]$. Since $v\in V_1$ and $|D(v)|=deg_T(v)-1=p-1$, $T[D[v]]=K_{1,p-1}$ with center $v$, and so $T=K_{1,p-1}\uplus_{vw}T'$.

We first show that $r_p(T')=p+1$. Due to $w\in V_2\subseteq \ma{D}_T$, $\ma{D}_T\cap V(T')$ is a $p$-dominating
set of $T'$, and so
 \begin{equation}\label{e4.10}
 \g_p(T)=|\ma{D}_T|=|\ma{D}_T\cap V(T')|+(p-1)\geq \g_p(T')+p-1.
 \end{equation}
Let $X'$ be an $\eta_p(T')$-set and $X=X'\cup D(v)$. Then, by (\ref{e4.10}),
 $$
 |X|=|X'|+|D(v)|=(\g_p(T')-1)+(p-1)<\g_p(T).
 $$
If $w\in X'$, then $r_p(T')=\eta_p(T')=\eta_p(V(T'),X',T')=
\eta_p(V(T),X,T)\geq r_p(T)=p+1$ by (\ref{e2.1})$\sim$(\ref{e2.3}) and Theorem \ref{thm2.2}. If $w\notin X'$, then $\eta_p(V(T),X,T)\geq p+2$ by $\mu_p(w,\ma{D}_T,T)\geq
p+2$ and Theorem~\ref{thm3.5}, and so
 $
 r_p(T')=\eta_p(V(T'),X',T')=\eta_p(V(T),X,T)-1\geq p+1.
 $
 We know $r_p(T')=p+1$ from Theorem~\ref{thm1.2}.

By Theorem~\ref{thm3.2}, $r_p(T')=p+1$ implies that $T'$ has a
unique $\g_p(T')$-set $\ma{D}_{T'}$.
We claim that $\ma{D}_{T'}=\ma{D}_T\cap V(T')$. To be
contrary, then $|\ma{D}_T\cap V(T')|\geq\g_p(T')+1$ since
$\ma{D}_T\cap V(T')$ is a $p$-dominating set of $T'$, and so, by (\ref{e4.10}),
 \begin{equation}\label{e4.11}
 |\ma{D}_{T'}\cup D[v]|=\g_p(T')+p\leq |\ma{D}_{T}\cap V(T')|-1+p=|\ma{D}_T|.
 \end{equation}
Since $\ma{D}_{T'}\cup D[v]$ is a
$p$-dominating set of $T$, (\ref{e4.11}) implies that $\ma{D}_{T'}\cup D[v]$ is a $\g_p(T)$-set
different to $\ma{D}_T$, a contradiction. The claim holds.

It is easily seen that $\ell_p(T')=\ell_p(T)-1$ from $\ma{D}_{T'}=\ma{D}_T\cap V(T')$ and $v\in PN_p(w,\ma{D}_T,T)$. Applying the induction on $T'$, $T'\in \mathscr{T}_p$ and
$V(T')\cap \ma{D}_T=\ma{D}_{T'}=A(T')$ by Lemma \ref{lem4.2}. Due to $w\in
\ma{D}_T\cap V(T')=A(T')$,  $T$ can be obtained from $T'$ by $\mathscr{O}_1$, and so $T\in \mathscr{T}_p$.

\textbf{Case 2.} $\mu_p(w,\ma{D}_T,T)=p+1$.

Since $r_p(T)=p+1$ and $w\in \ma{D}_T$, $w$ is not a stem of $T$, and so $C(w)\subseteq V_1$. Note that $N_T(w)=C(w)\cup \{x\}$ and
$\mu_p(w,\ma{D}_T,T)=p+1$. We can conclude from (\ref{e2.4}) that either $|C(w)|=2$ and
$x\in \ma{D}_T$ or $|C(w)|=1$ and $x\notin \ma{D}_T \cup
PN_p(w,\ma{D}_T,T)$ . We distinguish the following two subcases.

 \textbf{Subcase 2.1.} $|C(w)|=2$ and $x\in \ma{D}_T$.

Let $T'=T-D[w]$. Since $C(w)\subseteq V_1$ and $|C(w)=2|$, $T[D[w]]=F_{p-1}$ with center $w$, and so $T=F_{p-1}\uplus_{wx}T'$.

Due to $x\in \ma{D}_T$,
$x\notin PN_p(w,\ma{D}_T,T)$, and so, by \textbf{(ii)} of Theorem~\ref{thm3.3}, we have $r_p(T')=p+1$ and
$\ma{D}_{T'}=\ma{D}_T\cap V(T')$. So $\ell_p(T')=\ell_p(T)-|C(w)|<\ell_p(T)$.
Applying
the induction on $T'$, $T'\in \mathscr{T}_p$ and, by Lemma \ref{lem4.2},
$A(T')=\ma{D}_{T'}=\ma{D}_T\cap V(T')$.

Since $x\in \ma{D}_T$ and $r_p(T)=p+1$, by (\ref{e2.4})$\sim$(\ref{e2.6}) and Theorem \ref{thm2.4}, we have
 \begin{eqnarray*}
 |PN_p(x,\ma{D}_T,T)|+\max\{0,p-|N_T(x)\cap \ma{D}_T|\}&=&\mu_p(x,\ma{D}_T,T)\\
 &\geq&\mu_p(T)\geq r_p(T)=p+1,
 \end{eqnarray*}
which implies $|PN_p(x,\ma{D}_T,T)|\geq \min\{p+1,|N_T(x)\cap
\ma{D}_T|+1\}$. From $w\in N_T(x)\cap \ma{D}_T$ and
$A(T')=\ma{D}_T\cap V(T')$, we know that
 $$
 \mbox{$|PN_p(x,A(T'),T')|=|PN_p(x,\ma{D}_{T},T)|$\ \ and\ \ $|N_{T'}(x)\cap A(T')|=|N_T(x)\cap \ma{D}_T|-1$.}
 $$
Hence $|PN_p(x,A(T'),T')|\geq \min\{p+1,|N_{T'}(x)\cap A(T')|+2\}$.
Thus $T$ can be obtained from $T'$ by $\mathscr{O}_3$ and so
$T\in \mathscr{T}_p$.

\textbf{Subcase 2.2.} $|C(w)|=1$ and $x\notin \ma{D}_T \cup
PN_p(w,\ma{D}_T,T)$.

Let $T'=T-D[x]$.
By $x\in V_3$, $C(x)\subseteq V_0\cup V_1\cup V_2$. Note that the
unique neighbor of each vertex in $V_0$ is a $p$-private neighbor
with respect to $\ma{D}_T$ (this means that $\ma{D}_T\cap V_1=\emptyset$) and, to $p$-dominate $V_1$, all neighbors of each vertex in
$V_1$ must belong to $\ma{D}_T$ by (\ref{e4.9}). So we can derive from $x\notin \ma{D}_T \cup PN_p(w,\ma{D}_T,T)$ that $C(x)\subseteq V_2$ and $|C(x)|\geq p$.
By the
choice of $P=xwvu$ and $|C(w)|=1$, all vertices in $C(x)$ have degree $2$
in $T$. Hence $T[D[x]]=F_{t,p-1}$ with center $x$, where $t=|C(x)|\geq
p$. By $deg_T(x)\geq p>1=deg_T(r)$, we know that $x\neq r$ and so $x\in D(r)$. Let $y$ be the father of $x$, then $T=F_{t,p-1}\uplus_{xy}T'$.

Since $x\notin \ma{D}_T \cup
PN_p(w,\ma{D}_T,T)$, $\ma{D}_T\cap V(F_{t,p-1})$ (resp. $\ma{D}_T\cap V(T')$) is a $p$-dominating set of
$F_{t,p-1}$ (resp. $T'$). So,
 $$
 \g_p(T)=|\ma{D}_T|=|\ma{D}_T\cap V(F_{t,p-1})|+|\ma{D}_T\cap V(T')|\geq \g_p(F_{t,p-1})+\g_p(T').
 $$
Further, $\g_p(T)= \g_p(F_{t,p-1})+\g_p(T')$ since the union
between a $\g_p(F_{t,p-1})$-set and a $\g_p(T')$-set is a $p$-dominating
set of $T$. Hence $\ma{D}_T\cap V(F_{t,p-1})$ (resp. $\ma{D}_T\cap V(T')$) is a $\g_p(F_{t,p-1})$-set (resp. $\g_p(T')$-set).

If $\g_p(T')\leq p$, then we can check easily that $T'$ is a star
$K_{1,p}$ from $x\notin \ma{D}_T \cup
PN_p(w,\ma{D}_T,T)$ and $r_p(T)=p+1$. So $T$ can be
obtained from $K_{1,p}$ by $\mathscr{O}_4$ and $T\in
\mathscr{T}_p$.

If $\g_p(T')\geq p+1$, then, by Corollary \ref{cor2.3}, $p+1=r_p(T)\leq
r_p(T')$. It follows that $r_p(T')=p+1$ by Theorem~\ref{thm1.2}, and so
$\ma{D}_{T'}=\ma{D}_T\cap V(T')$ by Theorem \ref{thm3.2}. Note that $\ell_p(T')=\ell_p(T)-\ell_p(F_{t,p-1})< \ell_p(T)$. Applying the induction on $T'$, $T'\in
\mathscr{T}_p$ and $A(T')=\ma{D}_{T'}=\ma{D}_T\cap V(T')$ by Lemma \ref{lem4.2} and Theorem \ref{thm3.2}. Hence
$T$ can be obtained from $T'$ by $\mathscr{O}_4$ and $T\in
\mathscr{T}_p$.
\end{pf}

 Lemmas \ref{lem4.2} and \ref{lem4.3} imply that Theorem \ref{thm4.4} is true.

\begin{thm}\label{thm4.4}
For an integer $p\geq 3$ and a tree $T$, $r_p(T)=p+1$ if and only if
$T\in \mathscr{T}_p$.
\end{thm}



\end{document}